\RequirePackage{ifpdf}
\ifpdf 
\documentclass[pdftex]{sigma}
\else
\documentclass{sigma}
\fi

\newcommand{\sgn}{\operatorname{sgn}}
\newcommand{\pfaff}{\mathop{\mathrm{pfaf\/f}}}
\newcommand{\id}{\operatorname{id}}

\begin{document}
\allowdisplaybreaks
\renewcommand{\PaperNumber}{085}

\FirstPageHeading

\renewcommand{\thefootnote}{$\star$}

\ShortArticleName{Multivariable Christof\/fel--Darboux Kernels}

\ArticleName{Multivariable Christof\/fel--Darboux Kernels\\ and Characteristic Polynomials \\
of Random Hermitian Matrices\footnote{This paper is a contribution
to the Vadim Kuznetsov Memorial Issue ``Integrable Systems and
Related Topics''. The full collection is available at
\href{http://www.emis.de/journals/SIGMA/kuznetsov.html}{http://www.emis.de/journals/SIGMA/kuznetsov.html}}}

\Author{Hjalmar ROSENGREN}

\AuthorNameForHeading{H. Rosengren}

\Address{Department of Mathematical Sciences, Chalmers University
of Technology\\ and G\"oteborg
 University, SE-412~96 G\"oteborg, Sweden}

 \Email{\href{mailto:hjalmar@math.chalmers.se}{hjalmar@math.chalmers.se}}

\URLaddress{\href{http://www.math.chalmers.se/~hjalmar/}{http://www.math.chalmers.se/\~{}hjalmar/}}

\ArticleDates{Received October 11, 2006; Published online December 04, 2006}

\Abstract{We study multivariable Christof\/fel--Darboux kernels,
which may be viewed as reproducing kernels for antisymmetric
orthogonal polynomials, and also as
 correlation functions for products of characteristic polynomials of random Hermitian matrices.
 Using their interpretation as reproducing kernels, we obtain simple proofs of Pfaf\/f\/ian
  and determinant formulas, as well as Schur polynomial expansions, for such kernels.
In subsequent work, these results are applied in combinatorics
(enumeration of marked shifted tableaux) and number theory
(representation of integers as sums of squares).}


\Classification{15A15; 15A52; 42C05}

\begin{flushright}
\it Dedicated to the memory of Vadim Kuznetsov
\end{flushright}

\section{Introduction}

The Christof\/fel--Darboux kernel plays an important role in the
theory of one-variable orthogonal polynomials. In the present
work, we study a multivariable extension, which can be viewed as a
reproducing kernel for \emph{anti-symmetric}  polynomials. As is
explained below, our original motivation came from a very special
case, having applications in
 number theory (sums of squares) and combinatorics (tableaux enumeration).
More  generally, this kind of  kernels occur in random matrix
theory as correlation functions for products of characteristic
polynomials of random Hermitian matrices. The purpose of the
present note is to highlight a number of  useful identities for
such kernels. Although, as we will
 make clear, the main results can be found in the literature,
they are scattered in work belonging to dif\/ferent disciplines,
so it seems useful to collect them in one place. Moreover, our
proofs, with the interpretation as reproducing kernels, are new
and conceptually very simple.

We f\/irst recall something of the one-variable theory.  Let
\[
f\mapsto \int f(x)\,d\mu(x)
\]
be a linear functional def\/ined on polynomials of one variable.
We denote by  $V_n$ the  space  of polynomials of degree at most
$n-1$. Assuming that
 the pairing
\begin{displaymath}
\langle f,g\rangle = \int f(x)g(x)\,d\mu(x)
\end{displaymath}
is non-degenerate on each  $V_n$ (for instance, if it is positive
def\/inite), there exists a corresponding system
 $(p_k(x))_{k=0}^\infty$  of monic orthogonal polynomials.
We may then introduce the Christof\/fel--Darboux kernel
\begin{displaymath}
K(x,y)=\sum_{k=0}^{n-1}\frac{p_k(x)p_k(y)}{\langle
p_k,p_k\rangle},
\end{displaymath}
which is the \emph{reproducing kernel} of $V_n$,  that is, the
unique function  such that $(y\mapsto K(x,y))\in V_n$ and
\begin{displaymath}
f(x)=\int f(y)K(x,y)\,d\mu(y),\qquad f\in V_n.
\end{displaymath}
The Christof\/fel--Darboux formula  states that
\begin{displaymath}
K(x,y)=\frac 1{\langle
p_{n-1},p_{n-1}\rangle}\frac{p_n(x)p_{n-1}(y)-p_{n-1}(x)p_n(y)}{x-y}.
\end{displaymath}

More generally, let $V_n^m$, $0\leq m\leq n$, denote the $m$th
exterior power of $V_n$. It will be identif\/ied with the space of
antisymmetric polynomials $f(x)=f(x_1,\dots,x_m)$ that are of
degree at most $n-1$ in each $x_i$. Writing
\begin{displaymath}
\int f(x)\,d\mu_m(x)=\frac 1{m!}\int
f(x_1,\dots,x_m)\,d\mu(x_1)\dotsm d\mu(x_m),
\end{displaymath}
we equip $V_n^m$   with the pairing
\begin{displaymath}
\langle f,g\rangle_{V_n^m}=\int f(x)g(x)\,d\mu_m(x).
\end{displaymath}
Equivalently, in terms of the spanning vectors $\det\limits_{1\leq
i,j\leq m}(f_j(x_i))$, $f_j\in V_n$,
\begin{equation}\label{dsp}
\left\langle \det_{1\leq i,j\leq m}(f_j(x_i)),\det_{1\leq i,j\leq
m}(g_j(x_i)) \right\rangle_{ V_n^{ m}}=\det_{1\leq i,j\leq
m}(\langle f_i,g_j\rangle_{V_n}).
\end{equation}

Every element of $V_n^m$ is divisible by the polynomial
\begin{displaymath}\Delta(x)=\prod_{1\leq i<j\leq m}(x_j-x_i).
\end{displaymath}
The map $f\mapsto f/\Delta$ is an isometry from $V_n^m$ to the
space $W_n^m$, consisting of \emph{symmetric} polynomials in
$x_1,\dots,x_m$ that are of degree at most $n-m$ in each $x_i$,
equipped with the pairing
\begin{displaymath}
\langle f,g\rangle_{W_n^m}=\frac 1{m!}\int
f(x)g(x)\Delta(x)^2\,d\mu(x_1)\dotsm d\mu(x_m).
\end{displaymath}
We remark that, normalizing
\begin{equation}\label{pm}
\Delta(x)^2\,d\mu(x_1)\dotsm d\mu(x_m)
\end{equation}
to a probability distribution, it def\/ines an \emph{orthogonal
polynomial ensemble} (the term is sometimes used also for the more
general weights $|\Delta(x)|^\beta$). Such ensembles are important
in a variety of contexts, including the theory of random Hermitian
matrices \cite{ko}.

An orthogonal basis of $V_n^m$ is given by
$(p_S(x))_{S\subseteq[n],\,|S|=m}$, where
\begin{displaymath}
p_S(x)=\det_{1\leq i\leq m,\,j\in S}(p_{j-1}(x_i)).
\end{displaymath}
Here and throughout, we write $[n]=\{1,\dots,n\}$, and we assume
that the
 columns are ordered in  the natural way.
Indeed, it follows from \eqref{dsp} that
\begin{equation}\label{po}
\langle p_S,p_T\rangle=\delta_{ST}\prod_{i\in S}\langle
p_{i-1},p_{i-1}\rangle.
\end{equation}

We  denote by $\Delta(x)\Delta(y)K_m(x,y)$ the reproducing kernel
of $V_n^m$, that is,  the unique element of $V_n^m\otimes V_n^m$
such that
\begin{equation}\label{rp}
f(x)=\int f(y) \Delta(x)\Delta(y) K_m(x,y)\,d\mu_m(y),\qquad f\in
V_n^m.
\end{equation}
Equivalently, $K_m(x,y)$ is the reproducing kernel of $W_n^m$.

It is easy to see that
\begin{equation}\label{kw}
K_m(x,y)=\frac{1}{\Delta(x)\Delta(y)}\sum_{S\subseteq [n],\,
|S|=m}\frac{p_S(x)p_S(y)}{\langle p_S,p_S\rangle}
=\frac{\det\limits_{1\leq i,j\leq
m}(K(x_i,y_j))}{\Delta(x)\Delta(y)}.
\end{equation}
Indeed, it  follows from \eqref{dsp} and \eqref{po} that both
sides satisfy \eqref{rp}. The equality of the two expressions can
also be derived
 using the Cauchy--Binet formula.

The following elegant integral formula was recently obtained by
Strahov and Fyodorov \cite{sf}. In Section~\ref{s2} we will give a
simple proof, using the interpretation of the left-hand side as a
reproducing kernel.

\begin{proposition}[Strahov and Fyodorov]\label{il}
One has
\begin{equation}\label{ki}K_m(x,y)=\frac{1}{\prod\limits_{i=1}^n\langle p_{i-1},p_{i-1}\rangle}\int
\prod_{\substack{1\leq j\leq 2m\\[0.5mm] 1\leq k\leq n-m}}(z_j-w_k)\,
\Delta(w)^2\,d\mu_{n-m}(w), \end{equation} where
\begin{equation}\label{z}
(z_1,\dots,z_{2m})=(x_1,\dots,x_m,y_1,\dots,y_m).
\end{equation}
\end{proposition}

As a consequence, the
 obvious $S_m\times S_m\times \mathbb Z_2$ symmetry of $K_m(x,y)$ extends to a non-trivial $S_{2m}$ symmetry.

\begin{corollary}\label{c}
The polynomial $K_m(z)=K_m(z_1,\dots,z_{2m})$ is symmetric in all
its variables.
\end{corollary}

The motivation for the work of Strahov and Fyodorov is an
interpretation of the integral~\eqref{ki}  as a correlation
function for the product of characteristic polynomials of random
Hermitian matrices. The case $x=y$ is of particular interest,
since \eqref{ki} may then be written as
\begin{gather}\Delta(x_1,\dots,x_m)^2K_m(x_1,\dots,x_m,x_1,\dots,x_m)\nonumber\\
\qquad{}=\frac{1}{(n-m)!\prod\limits_{i=1}^n\langle
p_{i-1},p_{i-1}\rangle}
\int\Delta(x_1,\dots,x_n)^2\,d\mu(x_{m+1})\dotsm d\mu(x_{n}),
\label{kcf}
\end{gather}
which exhibits  $\Delta(x)^2K_m(x,x)$  as a correlation function
for the measure \eqref{pm}. In that case, the determinant formula
\eqref{kw} is classical, see   \cite[Chapter 5]{d}. In our
opinion, our proof of Proposition~\ref{il} is more illuminating
than the inductive proof usually given for the special case $x=y$.
For applications of Proposition \ref{il} and for further related
results, see \cite{bds,bs,bh}.

An alternative determinant formula for the integral  \eqref{ki}
may be obtained by combining  two classical results of
Christof\/fel \cite[Theorem 2.7.1]{is}  and Heine \cite[Theorem
2.1.2]{is}, see also \cite{bh}.
 We will obtain it as a by-product of the proof of Proposition \ref{il}.

\begin{proposition}\label{th1}
In the notation above,
\begin{displaymath}
K_m(z)=\frac{1}{\prod\limits_{i=1}^m\langle
p_{n-i},p_{n-i}\rangle} \frac{\det\limits_{1\leq i,j\leq
2m}(p_{n-m+j-1}(z_i))}{\Delta(z)}.
\end{displaymath}
\end{proposition}

More generally, such a  determinant formula holds for the
integrals
\begin{displaymath}\int
\prod_{\substack{1\leq j\leq m\\[0.5mm] 1\leq k\leq n}}(z_j-w_k)\,
\Delta(w)^2\,d\mu_{n}(w),\end{displaymath} but we focus  on the
case when $m$ is even.

The fact that the ``two-point'' determinant in \eqref{kw} and the
``one-point'' determinant in  Proposition \ref{th1} agree can also
be derived from the work of Lascoux \cite[Propositions 8.4.1 and
8.4.3]{l}, see \cite[Exercise 8.33]{l} for the integral formula in
this context. Note that the reproducing pro\-per\-ty  mentioned by
Lascoux, and given explicitly in  \cite[Proposition 3]{lh}, is of
a dif\/ferent nature from~\eqref{rp},
 pertaining to integration against \emph{one-variable} polynomials.

The special case of Proposition \ref{th1} obtained by subtracting
the $k$th row from the $(m+k)$th, for $1\leq k\leq m$, and then
letting $z_{m+k}\rightarrow z_k$, gives the following formula for
the correlation function~\eqref{kcf}.

\begin{corollary}
In the notation above,
\begin{displaymath}K_m(x,x)
= \frac{(-1)^{\frac 12m(m-1)}}{\prod\limits_{i=1}^m\langle
p_{n-i},p_{n-i}\rangle\Delta(x)^4}\,\det_{1\leq i,j\leq
2m}\left(\begin{cases}p_{n-m+j-1}(x_i), & 1\leq i\leq m,
\\ p_{n-m+j-1}'(x_i), & m+1\leq i\leq 2m
\end{cases}\right).
\end{displaymath}
\end{corollary}

Next, we give Pfaf\/f\/ian
 formulas for $K_m$. As is explained below, they can be deduced from results of  Ishikawa and Wakayama \cite{iw2}, Lascoux \cite{l2}, and Okada \cite{o}. Nevertheless, we will give an independent proof, using Corollary~\ref{c}.
Recall that the \emph{Pfaffian} of a skew-symmetric
even-dimensional matrix is given by
\begin{displaymath}\pfaff_{1\leq i,j\leq 2m}(a_{ij})
=\frac{1}{2^mm!}\sum_{\sigma\in S_{2m}}\sgn(\sigma)\prod_{i=1}^m
a_{\sigma(2i-1),\sigma(2i)}.\end{displaymath}

\begin{proposition}\label{th2}
For any choice of square roots  $\sqrt {z_i}$,
\begin{displaymath}
K_m(z)=\frac 1{\prod\limits_{1\leq i<j\leq 2m}(\sqrt {z_j}-\sqrt
{z_i})} \pfaff_{1\leq i,j\leq 2m}\left((\sqrt {z_j}-\sqrt
{z_i})K(z_i,z_j)\right).
\end{displaymath}
Moreover, for any choice of $\zeta_i$ such that
\begin{gather}\label{zeta}
\zeta_i+\zeta_i^{-1}=z_i+2,
\\
K_m(z)=\frac {\prod\limits_{i=1}^{2m}
\zeta_i^{m-1}}{\prod\limits_{1\leq i<j\leq 2m} ( {\zeta_j}-
{\zeta_i})}\pfaff_{1\leq i,j\leq 2m}\left(( {\zeta_j}-
{\zeta_i})K(z_i,z_j)\right). \nonumber
\end{gather}
\end{proposition}

Note that \eqref{zeta} implies
\begin{equation}\label{bb}z_j-z_i=-\frac 1{{\zeta_i\zeta_j}}\,({\zeta_j}-{\zeta_i})(1-{\zeta_i\zeta_j}).\end{equation}

In the special case $z=(x,x)$, choosing
\begin{gather*}(\sqrt{z_1},\dots,\sqrt{z_{2m}})=(-\sqrt{x_1},\dots,-\sqrt{x_m},\sqrt{x_1},\dots,\sqrt{x_m}),
\\
(\zeta_1,\dots,\zeta_{2m})=(\xi_1^{-1},\dots,\xi_m^{-1},\xi_1,\dots,\xi_m),
 \end{gather*}
Proposition \ref{th2} reduces to special cases of the  identity
\begin{displaymath}
\pfaff\left(\begin{matrix}A & B\\ -B &
-A\end{matrix}\right)=(-1)^{\frac12 m(m-1)}\det(B-A).
\end{displaymath}
The general case seems to lie deeper.

Proposition \ref{th2} is actually equivalent to  Proposition
\ref{p} below, which can be deduced from known  results. Indeed,
rewriting \eqref{dpi} using \cite[Theorem 4.2]{o} and \eqref{dpia}
using the case $n=m$ of \cite[Theorem 4.7]{o}, we can easily see
the  resulting expressions to agree. More explicitly, this
identity appears in \cite{l2}, with a simple proof. The
Pfaf\/f\/ian \eqref{dpib} can be treated similarly, or else shown
to agree with \eqref{dpia} by means of a result of Ishikawa and
Wakayama \cite[Theorem 5.1]{iw2}, see Remark~\ref{rem} below.

  We remark that  the relevant results of \cite{iw2},
\cite{l2} and \cite{o} are  closely related to Sundquist's
identities \cite{su}, see also \cite{iotz} and \cite{iw3}.
Moreover,
  \eqref{dpi} is related to the  Izergin--Korepin
  determinant for the partition function of the six-vertex model \cite{i},
  which, as well as the  Pfaf\/f\/ians in \eqref{dpio}, has  applications to alternating sign matrices \cite{k1,k2,o2,ze}.

\begin{proposition}\label{p}
Let $a_i$ and $z_i$, $1\leq i\leq 2m$, be free variables, and let
$\zeta_i$ be as in \eqref{zeta}. Moreover, let $S\subseteq [2m] $
be an arbitrary subset of cardinality $m$. Then,
\begin{subequations}\label{dpio}
\begin{gather}\label{dpi}(-1)^{\frac 12m(m+1)+\sum\limits_{s\in S}s}\prod_{i\in S,\,j\notin S}(z_j-z_i)
\det_{i\in S,\,j\notin S}\left(\frac{a_j-a_i}{z_j-z_i}\right)\\
\label{dpia}\qquad{}=
\prod_{1\leq i<j\leq 2m}(\sqrt{z_i}+\sqrt{z_j})\pfaff_{1\leq i,j\leq 2m}\left(\frac{a_j-a_i}{\sqrt{z_j}+\sqrt{z_i}}\right)\\
\label{dpib}\qquad{}=\prod_{i=1}^{2m}\zeta_i^{1-m}\prod_{1\leq
i<j\leq 2m}(1-\zeta_i\zeta_j) \pfaff_{1\leq i,j\leq
2m}\left(\frac{a_j-a_i}{1-\zeta_i\zeta_j}\right).
\end{gather}
\end{subequations}
\end{proposition}

Note that Proposition \ref{th2} is a special case of Proposition
\ref{p} when $a_i=p_n(z_i)/p_{n-1}(z_i)$,
$\{x_1,\dots,x_m\}=\{z_i\}_{i\in S}$ and
$\{y_1,\dots,y_m\}=\{z_i\}_{i\notin S}$.  Conversely, it is not
hard to deduce  Proposition \ref{p} from  Proposition \ref{th2}
using an interpolation argument.

\begin{remark}
Strahov and Fyodorov  found one-point and two-point determinant
formulas for more general
   correlation functions than those of Proposition \ref{il}. In the case of
\begin{displaymath}\int
\prod_{\substack{1\leq j\leq m\\[0.5mm] 1\leq k\leq n}}\frac{x_j-z_k}{y_j-z_k}\,
\Delta(z)^2\,d\mu_{n}(z),\end{displaymath} their two-point formula
is of the  type  \eqref{dpi}, see  \cite[Proposition 4.2]{sf}.
Thus, Proposition \ref{p} implies Pfaf\/f\/ian formulas for this
correlation function.
\end{remark}

In \cite{r3}, we  need an elementary result on the expansion of
the kernel $K_m(x,y)$ into Schur polynomials
$s_\lambda(x)s_\mu(y)$. Since we have not found a suitable
reference, we include it here. The proof is given in Section
\ref{ss}.

\begin{proposition}\label{sp} One has
\begin{displaymath}
K_m(x,y)=\sum_{\substack{0\leq \lambda_m\leq\dots\leq \lambda_1\leq n-m\\[0.5mm] 0\leq \mu_m\leq\dots\leq \mu_1\leq n-m}}
\frac{(-1)^{\sum\limits_{i=1}^m(\lambda_i+\mu_i)}}{\prod\limits_{i=1}^n\langle
p_{i-1},p_{i-1}\rangle} \det_{i\in [n]\setminus S,j\in
[n]\setminus T}(c_{i+j-2})\,s_\lambda(x)s_\mu(y),
\end{displaymath} where
\begin{displaymath}S=\{\lambda_k+m+1-k;\,1\leq k\leq m\},\qquad
T=\{\mu_k+m+1-k;\,1\leq k\leq m\},\end{displaymath} and
\begin{displaymath}
c_k=\int x^k\,d\mu(x).
\end{displaymath}
\end{proposition}

Finally, let us describe our original motivation, which comes from
applications that seem completely unrelated to random matrix
theory. Motivated by the theory of af\/f\/ine super\-algebras, Kac
and Wakimoto \cite{kw} conjectured certain  new formulas for the
number of representations of an integer as the sum of $4m^2$ or
$4m(m+1)$ triangular numbers. These conjectures were f\/irst
proved by Milne \cite{mi1,mi3,mi2}, and later by Zagier \cite{z}.
In \cite{ro}, we re-derived and generalized the Kac--Wakimoto
identities using elliptic Pfaf\/f\/ian evaluations. Extension of
this analysis from triangles to squares leads to formulas
involving \emph{Schur $Q$-polynomials} \cite{ma} evaluated at the
point $(1,\dots,1)$. (More precisely, these polynomials are
normally labelled by positive integer partitions. Here, we need an
extension  to the case when some indices are \emph{negative}.)
Later, we realized that the resulting sums of squares formulas are
equivalent to those of Milne \cite{mi2}. Seeing this is far from
obvious and requires an identif\/ication
  of the Schur $Q$-polynomials with  kernels $K_m(z)$,
  where the underlying  polynomials $p_k(x)$ are continuous dual Hahn polynomials.
  The key fact for obtaining this identif\/ication is the second part of Proposition~\ref{th2}.
  We refer to \cite{r2} for  applications of the results above
  to Schur $Q$-polynomials and marked shifted tableaux, and to \cite{r3} for the relation to sums of squares.

\section{Proof of Propositions \ref{il} and \ref{th1}}\label{s2}

\begin{lemma}\label{l}
Let $\phi:\,V_n^m\rightarrow V_n^{n-m}$ be def\/ined by
\begin{displaymath}(\phi f)(x)=\int f(y)\Delta(y,x)\,d\mu_m(y). \end{displaymath}
Then,
\begin{displaymath}(\phi p_S)(x)=(-1)^{\frac12m(m+1)+\sum\limits_{s\in S}s}
\prod_{i\in S}\langle p_{i-1},p_{i-1}\rangle\cdot p_{S^c}(x),
\end{displaymath} where $S^c=\{ 1,\dots,n\}\setminus S$.
\end{lemma}

A corresponding statement holds  when
 $V_n$ is a general $n$-dimensional vector space and
$\Delta$  an element of the one-dimensional space $(V_n^\ast)^n$
\cite[\S~8.5]{b}, the map
 $\phi:\,V_n^m\rightarrow (V_n^\ast)^{n-m}$ often being called \emph{Hodge star} or \emph{Poincar\'e isomorphism}.
 For completeness, we provide a proof in the present setting.

\begin{proof}
The Vandermonde determinant evaluation
\begin{equation}\label{vd}\Delta(x)=\det_{1\leq i,j\leq n}(p_{j-1}(x_i))\end{equation}
gives
\begin{gather*}(\phi p_S)(x)
=\sum_{\substack{\sigma:\,[m]\rightarrow S\\
\tau:\,[n]\rightarrow [n] }}\sgn(\sigma)\sgn(\tau)
\int\prod_{i=1}^m
p_{\sigma(i)-1}(y_i)p_{\tau(i)-1}(y_i)\,d\mu_m(y)\prod_{i=m+1}^{n}p_{\tau(i)-1}(x_i),
 \end{gather*}
where the sum is over bijections. By orthogonality, we may assume
\begin{displaymath}(\tau(1),\dots,\tau(n))=(\sigma(1),\dots,\sigma(m),\rho(1),\dots,\rho(n-m)), \end{displaymath}
with $\rho$  a bijection $[n-m]\rightarrow S^c$. One easily checks
that
\begin{displaymath}\sgn(\tau)=(-1)^{\frac12m(m+1)+\sum\limits_{s\in S}s}\sgn(\sigma)\sgn(\rho), \end{displaymath}
which gives indeed
\begin{gather*}(\phi p_S)(x)=(-1)^{\frac12m(m+1)
+\sum\limits_{s\in S}s}\frac{1}{m!}\sum_{\sigma:\,[m]\rightarrow
S}
\prod_{i=1}^m \langle p_{\sigma(i)-1},p_{\sigma(i)-1}\rangle\\
\phantom{(\phi p_S)(x)=}{} \times\sum_{\rho:\,[n-m]\rightarrow S^c}\sgn(\rho)\prod_{i=1}^{n-m}p_{\rho(i)-1}(x_i)\\
\phantom{(\phi p_S)(x)}{}=(-1)^{\frac12m(m+1)+\sum\limits_{s\in
S}s}\prod_{i\in S} \langle p_{i-1},p_{i-1}\rangle\,p_{S^c}(x).
\tag*{\qed}
\end{gather*}
\renewcommand{\qed}{}
\end{proof}

By iteration, it follows from Lemma \ref{l} that
\begin{displaymath}
\phi\circ\phi =(-1)^{m(n-m)}\prod_{i=1}^n\langle
p_{i-1},p_{i-1}\rangle\cdot\id_{V_n^m}.
\end{displaymath}
By means of $\Delta(w,x)=(-1)^{m(n-m)}\Delta(x,w)$, $x\in\mathbb
R^m$, $w\in\mathbb R^{n-m}$, this fact can be expressed as
\begin{displaymath}
f(x)=\frac{1}{\prod\limits_{i=1}^n\langle p_{i-1},p_{i-1}\rangle}
\int f(y)\Delta(x,w)\Delta(y,w)\,d\mu_m(y)d\mu_{n-m}(w),\qquad
f\in V_n^m.
 \end{displaymath}
By the uniqueness of the reproducing kernel, it follows that
\begin{equation}\label{idd}\Delta(x)\Delta(y)K_m(x,y)
=\frac{1}{\prod\limits_{i=1}^n\langle
p_{i-1},p_{i-1}\rangle}\int\Delta(x,w)\Delta(y,w)\,d\mu_{n-m}(w).
\end{equation} This is equivalent to Proposition \ref{il}.

\begin{remark}
The equation \eqref{idd}  can also be obtained as the special case
$l=n$ of the contraction formula
\begin{equation}\label{cf}\Delta(x)\Delta(y)K_m(x,y)
=\frac{(n-l)!(l-m)!}{(n-m)!}\int\Delta(x,w)\Delta(y,w)
K_l(x,w,y,w)\,d\mu_{l-m}(w),\end{equation} $0\leq m\leq l\leq n$.
Conversely, \eqref{cf} follows easily from~\eqref{idd}.
\end{remark}

To prove Proposition \ref{th1}, we note that
\begin{displaymath}
\Delta(x,y,w)=\frac{\Delta(x,y)\Delta(x,w)\Delta(y,w)}{\Delta(x)\Delta(y)\Delta(w)}.
\end{displaymath}
Applying this to \eqref{idd} gives
\begin{gather*}
K_m(x,y)=\frac{1}{\prod\limits_{i=1}^n\langle
p_{i-1},p_{i-1}\rangle} \frac{1}{\Delta(x,y)}
\int\Delta(x,y,w)\Delta(w)\,d\mu_{n-m}(w)\\
\phantom{K_m(x,y)}{}=\frac{1}{\prod\limits_{i=1}^n\langle
p_{i-1},p_{i-1}\rangle} \frac{(\phi\Delta)(z)}{\Delta(z)},
\end{gather*}
where, as in \eqref{z}, $z=(x,y)$. Next we observe that, by
\eqref{vd}, $\Delta(w)=p_S(w)$, with $S=[n-m]$. Lemma~\ref{l} then
gives indeed
\begin{displaymath}
K_m(x,y)=\frac{1}{\prod\limits_{i=1}^m\langle
p_{n-i},p_{n-i}\rangle} \frac{\det\limits_{1\leq i,j\leq
2m}(p_{n-m+j-1}(z_i))}{\Delta(z)}. \end{displaymath}

\section{Proof of Proposition  \ref{th2}}\label{s3}

Our main tool is the following elementary property of
Pfaf\/f\/ians, which we learned from an unpublished manuscript of
Eric Rains \cite{ra}.

\begin{lemma}[Rains]\label{rl}
For arbitrary $(a_{ij})_{1\leq i,j\leq 2m}$,
\begin{displaymath}\pfaff_{1\leq i,j\leq 2m}(a_{ij}-a_{ji})=\sum_{{S\subseteq [2m],\, |S|=m}}
(-1)^{\chi(S)}\det_{i\in S,j\notin S}(a_{ij}),
 \end{displaymath}
where $\chi(S)$ denotes the number of even elements in $S$.
\end{lemma}

For completeness, we sketch Rains' proof.

\begin{proof}
The left-hand side is given by
\begin{displaymath}\frac 1{2^m m!}\sum_{\sigma\in S_{2m}}\sgn(\sigma)\prod_{i=1}^{m}(a_{\sigma(2i-1),
{\sigma(2i)}}-a_{\sigma(2i),{\sigma(2i-1)}}) =\frac 1{
m!}\sum_{\sigma\in S_{2m}}\sgn(\sigma)\prod_{i=1}^{m}
a_{\sigma(2i-1),{\sigma(2i)}}.
\end{displaymath}
Write
\begin{displaymath}
\prod_{i=1}^{m} a_{\sigma(2i-1),\sigma(2i)}=\prod_{i\in
S}a_{i,\tau(i)},
\end{displaymath}
with $S=\{\sigma(1),\sigma(3),\dots,\sigma(2m-1)\}$ and $\tau$ a
bijection $S\rightarrow S^{c}$. Identifying $\tau$ as an element
of~$S_m$, using the natural orderings on  $S$ and $S^{c}$, it is
easy to check that $\sgn(\sigma)=(-1)^{\chi(S)}\sgn(\tau)$. Since
the map $\sigma\mapsto(S,\tau)$ is $m!$ to one,  we obtain indeed
\begin{gather*}\sum_{{S\subseteq [2m],\, |S|=m}}(-1)^{\chi(S)}\sum_{\tau}
\sgn(\tau)\prod_{i\in S}a_{i,\tau(i)}.\tag*{\qed}
 \end{gather*}
 \renewcommand{\qed}{}
\end{proof}

The following identity appeared as  \cite[Theorem A.1]{iw2}.

\begin{corollary}[Ishikawa and Wakayama]\label{iwp}
One has
\begin{gather*}
\pfaff_{1\leq i,j\leq 2m}\left(\frac{z_j-z_i}{a+b(x_i+x_j)+cx_ix_j}\right)\\
\qquad{}=(b^2-ac)^{\frac12m(m-1)}\prod_{1\leq i<j\leq 2m}\frac 1{a+b(x_i+x_j)+cx_ix_j}\\
\qquad\phantom{{}=}{}\times\sum_{{S\subseteq [2m],\, |S|=m}}
(-1)^{\chi(S)}\prod_{j\notin S}z_j\prod_{\substack{1\leq i<j\leq 2m\\[0.5mm] i,j\in S
\text{ \emph{or} }i,j\notin S}}(x_j-x_i)(a+b(x_i+x_j)+cx_ix_j).
\end{gather*}
\end{corollary}

This follows immediately from Lemma \ref{rl}, by use of the Cauchy
determinant
\begin{displaymath}\det_{1\leq i,j\leq m}\left(\frac{1}{a+b(x_i+y_j)+cx_iy_j}\right)
=(b^2-ac)^{\frac12m(m-1)}\frac{\prod\limits_{1\leq i<j\leq
m}(x_j-x_i)(y_j-y_i)}{\prod\limits_{i,j=1}^m(a+b(x_i+y_j)+cx_iy_j)},
\end{displaymath}
which is reduced to its more well-known special case $a=1$, $b=0$,
$c=-1$ through the elementary identity
  \begin{equation}\label{ei}
a+b(x_i+x_j)+cx_ix_j=\frac 1 c\left(ac-{b^2}+(c x_i+{b})( c
x_j+{b})\right).
\end{equation}
 The proof of Corollary \ref{iwp} given in \cite{iw2} is more complicated.

\begin{remark}\label{rem}
The  equality of \eqref{dpia} and \eqref{dpib} can be deduced by
applying Corollary \ref{iwp} to both Pfaf\/f\/ians,   using also
\eqref{bb}.
\end{remark}

We only need  Corollary \ref{iwp} in the case when $x_i=z_i$.
Then,
 the Pfaf\/f\/ian is given by
\begin{displaymath}(b^2-ac)^{m(m-1)}
\prod_{1\leq i<j\leq
2m}\frac{x_j-x_i}{a+b(x_i+x_j)+cx_ix_j}.\end{displaymath} Indeed,
one may use \eqref{ei} to reduce oneself  to the case $a=1$,
$b=0$, $c=-1$, which is the Pfaf\/f\/ian evaluation in
\cite[Proposition 2.3]{ste}.

\begin{corollary}\label{ssc}
One has
\begin{gather*}\sum_{\substack{S\subseteq [2m],\, |S|=m}}
(-1)^{\chi(S)}\prod_{j\notin S}x_j
\prod_{\substack{1\leq i<j\leq 2m\\[0.5mm] i,j\in S \text{ \emph{or} }i,j\notin S}}(x_j-x_i)(a+b(x_i+x_j)+cx_ix_j)\\
\qquad{}=(b^2-ac)^{\frac12m(m-1)} \prod_{1\leq i<j\leq
2m}(x_j-x_i).
 \end{gather*}
\end{corollary}

We are now ready to prove Proposition \ref{th2}. By Lemma
\ref{rl}, we have in general
\begin{displaymath}\pfaff_{1\leq i,j\leq 2m}((a_j-a_i)K(z_i,z_j))=\sum_{{S\subseteq [2m],\, |S|=m}}
(-1)^{\chi(S)}\prod_{j\notin S}a_j\det_{i\in S,j\notin
S}(K(z_i,z_j)), \end{displaymath} which, by Corollary \ref{c},
equals
\begin{equation}\label{gp}K_m(z)\sum_{{S\subseteq [2m],\, |S|=m}}
(-1)^{\chi(S)}\prod_{j\notin S}a_j\prod_{\substack{1\leq i<j\leq 2m\\[0.5mm]
 i,j\in S \text{{ or} }i,j\notin S}}(z_j-z_i). \end{equation}

Consider f\/irst the case $a_i=\sqrt z_i$. Then, by the case
$a=c=0$, $b=1$, $x_i= \sqrt {z_i}$ of  Corollary~\ref{ssc}, the
sum in \eqref{gp} equals
\begin{displaymath}\prod_{1\leq i<j\leq 2m}\left(\sqrt{z_j}-\sqrt{z_i}
\right). \end{displaymath}
 This yields the f\/irst part of Proposition \ref{th2}. Similarly, letting $a_i=\zeta_i$ and using \eqref{bb},
we can compute the sum in~\eqref{gp} by the case $a=1$, $b=0$,
$c=-1$, $x_i=\zeta_i$ of
 Corollary~\ref{ssc} as
\begin{displaymath}\prod_{i=1}^{2m}\zeta_i^{1-m}\prod_{1\leq i<j\leq 2m}\left(\zeta_j-\zeta_i\right).  \end{displaymath}
This completes the proof of Proposition \ref{th2}.

\section{Proof of Proposition \ref{sp}}
\label{ss}

When $(e_k)_{k=1}^n$ is a basis of $V_n$, let
$(e_S)_{S\subseteq[n],\, |S|=m}$ be the corresponding basis of
$V_n^m$ def\/ined by
\begin{displaymath}
e_S(x)=\det_{1\leq i\leq m,\,j\in S}(e_j(x_i)).
 \end{displaymath}
We then have the following general expansion formula. If
$e_i(x)=f_i(x)=p_{i-1}(x)$, this is \eqref{kw}.

\begin{lemma}\label{gel}
Let $(e_k)_{k=1}^n$ and  $(f_k)_{k=1}^n$ be arbitrary bases of
$V_n$. Then,
\begin{equation}\label{ge}
\Delta(x)\Delta(y)K_m(x,y)
=\sum_{\substack{S,T\subseteq[n]\\[0.5mm] |S|=|T|=m}}
(-1)^{\sum\limits_{s\in S}s+\sum\limits_{t\in T}t}\,
\frac{\det\limits_{i\in S^c,j\in T^c}(\langle e_i,f_j\rangle)}
{\det\limits_{1\leq i,j\leq n}(\langle
e_i,f_j\rangle)}\,e_S(x)f_T(y).
\end{equation}
\end{lemma}

\begin{proof}
It is enough to show that
\begin{displaymath}e_U(x)=\int e_U(y)R(x,y)d\mu_m(y),\qquad U\subseteq[n],\quad |U|=m, \end{displaymath}
where $R$ denotes the right-hand side of \eqref{ge}. Equivalently,
we need to show that
\begin{displaymath}\sum_{T\subseteq[n],\,|T|=m}
(-1)^{\sum\limits_{s\in S}s+\sum\limits_{t\in T}t}\,
\frac{\det\limits_{i\in S^c,j\in T^c} (\langle
e_i,f_j\rangle)\det\limits_{i\in U,j\in T}(\langle
e_i,f_j\rangle)} {\det\limits_{1\leq i,j\leq n}(\langle
e_i,f_j\rangle)}=\delta_{SU}.
\end{displaymath}
To verify this, write $U=\{u_1<\dots<u_m\}$,
$S^c=\{s_1'<\dots<s_{n-m}'\}$, and let $X$ be the ordered sequence
$(u_1,\dots,u_m,s_1',\dots,s_{n-m}')$. Consider the determinant
\begin{displaymath}D=\det_{i\in X,1\leq j\leq n}(\langle e_i,f_j\rangle). \end{displaymath}
On the one hand, reordering the rows gives
\begin{displaymath}D=(-1)^{\frac 12m(m+1)+\sum\limits_{s\in S}s}\delta_{SU}\det_{1\leq i,j\leq n}(\langle e_i,f_j\rangle).\end{displaymath}
 On the other hand, applying Laplace expansion to the f\/irst $m$ rows gives
\begin{displaymath}D=\sum_{T\subseteq[n],\,|T|=m}
(-1)^{\frac 12m(m+1)+\sum\limits_{t\in T}t}\det_{i\in U,j\in
T}(\langle e_i,f_j\rangle) \det_{i\in S^c,j\in T^c}(\langle
e_i,f_j\rangle). \end{displaymath} This completes the proof.
\end{proof}

 Consider the case of Lemma \ref{gel}
 when  $e_i(x)=f_i(x)=x^{i-1}$. Then,  $e_S(x)=\Delta(x)s_\lambda(x)$,
 where $S=\{s_1<\dots<s_m\}$ and $\lambda_i=s_{m+1-i}-i$. Thus,  noting also that
\begin{displaymath}
\det_{1\leq i,j\leq n}(\langle e_i,f_j\rangle) =\det_{1\leq
i,j\leq n}(\langle p_{i-1},p_{j-1}\rangle) =\prod_{i=1}^n\langle
p_{i-1},p_{i-1}\rangle, \end{displaymath} we obtain Proposition
\ref{sp}.

\subsection*{Acknowledgements}
I  thank Alain Lascoux and
 Eric Rains for  communicating  their unpublished works \cite{l2,ra},
as well as for several useful comments. The research  was
supported by the Swedish Science Research
 Council (Vetenskapsr\aa det).

\LastPageEnding

\end{document}